\documentclass[10pt]{article}
\usepackage{amsmath,amssymb,latexsym}
\usepackage{bbm} % bbold
\usepackage{epsfig}
\usepackage{lscape} 
\usepackage{array}

\def\bA{\boldsymbol{A}}

\def\bG{\boldsymbol{G}}
\def\bI{\boldsymbol{I}}
\def\bU{\boldsymbol{U}}

\def\bK{\boldsymbol{K}}

\def\bR{\boldsymbol{R}}

\def\bX{\boldsymbol{X}}
\def\bY{\boldsymbol{Y}}

\def\bPi{\boldsymbol{\Pi}}

\def\ba{\boldsymbol{a}}
\def\grasb{\boldsymbol{b}}
\def\bc{\boldsymbol{c}}

\def\bg{\boldsymbol{g}}
\def\bk{\boldsymbol{k}}
\def\bx{\boldsymbol{x}}
\def\by{\boldsymbol{y}}
\def\bz{\boldsymbol{z}}
\def\bu{\boldsymbol{u}}
\def\bv{\boldsymbol{v}}

\def\bt{\boldsymbol{t}}
\def\bs{\boldsymbol{s}}

\def\be{\boldsymbol{e}}

\def\bdelta{\boldsymbol{\delta}}
\def\bgamma{\boldsymbol{\gamma}}
\def\bzero{\boldsymbol{0}}
\def\bminfty{\boldsymbol{-\infty}}

\def\real{\mathbb{R}}

\def\integer{\mathbb{N}}

\def\a{\mathbbm{a}}
\def\c{\mathbbm{c}}

\def\mat1{\mathbb{I}}

\newcommand { \Iletter}[1] {I\kern-0.10em #1 }
\def\bit{\begin{itemize}}
\def\eit{\end{itemize}}
\def\ben{\begin{enumerate}}
\def\een{\end{enumerate}}
\def\bde{\begin{description}}
\def\ede{\end{description}}

\def\bar{\begin{array}}
\def\ear{\end{array}}
\def\beq{\begin{equation}}
\def\eeq{\end{equation}}
\def\bfi{\begin{figure}[hbt] \begin{center}}
\def\efi{\end{center} \end{figure}}
\def\noi{\noindent}
\def\bce{\begin{center}}
\def\ece{\end{center}}
\newcommand{\proof}{{\bf Proof. }}
\newcommand{\cqfd}{\hfill $\Box$}

\newtheorem {theo} {Theorem}[section]

\newtheorem {lemm} {Lemma}[section]

\newtheorem {defi} {Definition}[section]

\newtheorem {rem} {Remark} [section]

\newtheorem {resu} {Result} [section]

\begin{document}        

%\linenumbers
     
%\begin{frontmatter}                     

\title{Some Ideas to Test if a Polyhedron is Empty}
\author{Laurent Truffet \\
  IMT-A \\
Dpt. Automatique-Productique-Informatique \\  
Nantes, France \\
email: laurent.truffet@imt-atlantique.fr}

%\ead{chams.lahlou@imt-atlantique.fr}
%\author[emn,dap]{L. Truffet}
%\ead[url]{http://www.emn.fr/truffet}

%\cortext[cor1]{Corresponding author}
% \cortext[cor2]{Corresponding author}

%\address[emn]{Ecole des Mines de Nantes, La Chantrerie, 4 rue
%A. Kastler, BP 20722, Nantes 44307 Cedex 3, France}
%\address[dsee]{Dpt. Syst\`emes Energ\'etiques et Environnement}
%\address[dap]{Dpt. Automatique-Productique}

\maketitle 

\begin{abstract}
  In this paper we develop a pure algebraic method which provides an
  algorithm for testing emptyness of a polyhedron.
\end{abstract}

\noi

%\begin{keyword} 
%% keywords here, in the form: keyword \sep keyword

%% MSC codes here, in the form: \MSC code \sep code
%% or \MSC[2008] code \sep code (2000 is the default)

%\end{keyword}

%\end{frontmatter}

\section{Introduction} 
\label{secIntro}
Testing if a polyhedron is empty is a fundamental issue for linear
optimization theory. In the seminal work \cite{kn:Khachiyan79}
a geometric approach was proposed to solve this problem. Since
this work and to the best knowledge of the author only geometric
approaches were proposed (see e.g. \cite{kn:Karmarkar84}, \cite{kn:Chubanov15} among
many others). In this work we develop a method based only
on algebraic considerations. The main concepts deal with Moore-Penrose
inverse of a rectangular matrix (see Definition~\ref{defMP}) and some
results of arithmetic intervals \cite{kn:Moore79}. This method provides
an algorithm which seems to be new for testing emptyness of a
polyhedron.

\subsection{Main notations}

For any integer $k \geq 1$ we defined the set $[k]:=\{1, \ldots ,k\}$.

Bold letters represent matrices or column vectors. $^{t}(\cdot)$
denotes transpose operator. $\textsf{Mat}(\real,m,n)$ denotes the set
of all real $m \times n$-matrices.

We denote by $\bzero_{p,q} \in \textsf{Mat}(\real, p,q)$ the $p\times
q$-matrix whose elements are all zero. We denote by $\bI_{p} \in
\textsf{Mat}(\real,p,p)$ the $p\times p$-identity matrix. We denote
$\bzero$ (resp. $\bminfty$) the column vector whose components are all $0$
(resp. $-\infty$). The dimension is fixed by the context.

The natural order defined on $\real$ is denoted $\leq$. For any
integer $p \geq 2$ this total order is extended to the partial order
(called componentwise ordering), once again denoted $\leq$, defined on
$p$-dimensional vectors as follows.  If $\bx= \; ^{t}(x_{1}, \ldots,
x_{p})$ and $\by= \; ^{t}(y_{1}, \ldots, y_{p})$ then

\begin{equation}
\bx \leq \by \Leftrightarrow \forall i=1, \ldots, p: \; x_{i} \leq y_{i}.
  \end{equation}
The relation $\by \geq \bx$ means $\bx \leq \by$.

The vector $\be^{i}$ denotes the vector such that $\forall j, \;
e^{i}_{j}=1$ if $j=i$ and $0$ otherwise. Its dimension is determined by the context.

\subsection{Problem statement}
\label{sub-pb}

In this paper we consider the following set:
\begin{equation}
\mathcal{P}(\bA, \grasb):=\{ \bx \in \real^{n}: \bA \bx \leq \grasb\}.
\end{equation}
$\bA \in \textsf{Mat}(\real,m,n)$ and $\grasb \in \real^{m}$. \\
And we wonder wether this set is empty or not ? \\

\noi
We make the following assumptions (see Section~\ref{sec-Discuss} for a discussion): \\

\noi
\textrm{ASSUMPTION (A)}. Matrix $\bA$ has no null row. \\

\noi
\textrm{ASSUMPTION (B)}. We assume $m > n$ and that $\bA$ has full column-rank, that is
$\textsf{rk}(\bA)=n$. \\

\subsection{Basic concepts and results}
Our method is based on the following concepts, remarks and results.

\begin{defi}[Moore-Penrose inverse \cite{kn:Moore20}, \cite{kn:Penrose55}]
  \label{defMP}
  The Moore-Penrose inverse of the matrix $\bA \in \textsf{Mat}(\real,m,n)$ is the
  unique matrix $\bA^{+} \in \textsf{Mat}(\real,n,m)$ such that: \\
  \ben
\item $\bA \bA^{+} \bA = \bA$
\item $\bA^{+} \bA \bA^{+} = \bA^{+}$
\item $^{t}(\bA \bA^{+}) = \bA \bA^{+}$
  \item $^{t}(\bA^{+} \bA)= \bA^{+} \bA$
  \een

A matrix which satifies relations $1$ and $2$ will be
called a $\{1,2\}$-inverse.
  
  \end{defi}

We make the following remarks:

\bit
\item[R1] $\mathcal{P}(\bA, \grasb) \neq \emptyset \Leftrightarrow
\exists \bx, \exists \bc \leq \grasb, \bA \bx = \bc$ 

\item[R2] $\exists \bx, \bA \bx = \bc \Leftrightarrow \bA \bA^{+} \bc = \bc$, 
recalling that $\bA^{+} \in \textsf{Mat}(\real, n,m)$ denotes the 
Moore-Penrose inverse of $\bA$. 

\eit

Remark (R1) is obvious. Remark (R2) can be found in e.g. \cite{kn:James78}. As a main
consequence of (R1) and (R2) we have:
\begin{equation}
  \label{eq-nonempty}
\mathcal{P}(\bA, \grasb) \neq \emptyset \Leftrightarrow \mathfrak{A}: \exists \bc \leq \grasb,
\bA \bA^{+} \bc = \bc.
\end{equation}

Noticing that each component $c_{i}$ of $\bc$ belongs to the interval
$\mathcal{C}_{i}:=[-\infty, b_{i}]$, $i=1, \ldots, m$, we recall some results of interval
arithmetic \cite{kn:Moore79}.

Let $-\infty \leq s \leq t \leq +\infty$ and 
$-\infty \leq s' \leq t' \leq +\infty$ we define the intervals 
$[s,t]$ and $[s',t']$ as the following sets: 
$[s,t]:=\{a : s \leq a \leq t\}$ and $[s',t']:=\{a : s' \leq a \leq t'\}$. 

\bit
\item Interval addition. The addition of the intervals $[s,t]$ 
and $[s',t']$ is the new interval denoted by $[s,t] + [s',t']$ and 
defined by:
\begin{equation}
[s,t] + [s',t']:= \{a + a', \; a \in [s,t], \; a' \in [s',t']\} =[s+s',t+t'].
\end{equation}
Because the addition of reals is commutative and associative so is the 
interval addition. The interval $[0,0]$ being its neutral element.

\item Interval multiplication by a real. Let $z \in \real$ the 
multiplication of $z$ by the interval $[s,t]$ provides the new 
interval denoted $z \cdot [s,t]$ and defined by:
\begin{equation}
z \cdot [s,t] := \{z a, \; a \in [s,t]\} = \left\{ \bar{ll}
    [z s, z t] & \mbox{if $z > 0$} \\
      \mbox{$[ z t, z s]$} & \mbox{if $z < 0$} \\
      \mbox{$[0,0]$} & \mbox{if $z= 0$}
    \ear \right.
\end{equation}
\item Linear combination of two intervals. Let $z,z' \in \real$. The linear combination
  $z \cdot [s,t] + z' \cdot [s',t']$ is the set defined as:
  \begin{equation}
z \cdot [s,t] + z' \cdot [s',t']:= \{z a+ z'a', \; a \in [s,t], \; a' \in [s',t']\}.
    \end{equation}
\eit

Let $r \geq 1$, and let $I_{i}:=[s_i, t_i]$, $i=1, \ldots, r$ be a series of $r$
intervals. Let $\boldsymbol{\a}:= \; ^{t}(I_{1}, \ldots, I_{r})$ be the $r$-dimensional vector
of the intervals $I_i$, $i =1, \ldots,r$. Let us also define
$\bs:= \; ^{t}(s_{1}, \ldots, s_{r})$ and $\bt:= \;^{t}(t_{1}, \ldots, t_{r})$. We then have:
\begin{equation}
\boldsymbol{\a} =\{\ba : \bs \leq \ba \leq \bt \}.
\end{equation}

An interval $[s,t]$ is {\em thin} if $s=t$. By extension we say that
the interval vector $\boldsymbol{\a}$ is {\em thin} if $\bs=\bt$.

\begin{resu}[Beeck]
  \label{res-Beeck}
  For all $\bz \in \real^{r}$ we have:

  \begin{equation}
    \exists \ba \mbox{ s.t.} \; \bs \leq \ba \leq \bt \mbox{ and} \;
    ^{t}\bz \ba =0 \Leftrightarrow 0 \in \; ^{t}\bz \cdot \boldsymbol{\a}.
  \end{equation}
  Where
  \[
^{t}\bz \cdot \boldsymbol{\a}:=z_{1} \cdot I_{1}+ \ldots + z_{r} \cdot I_{r},
  \]
  with $z_{1} \cdot I_{1}+ \ldots + z_{r} \cdot I_{r}$ that generalizes the linear combination of
  two intervals as follows: $z_{1} \cdot I_{1}+ \ldots + z_{r} \cdot I_{r}:=\{\sum_{i=1}^{r} z_{i} a_{i}, \; a_{i} \in I_{i}, i=1, \ldots, r\}$.
  
\end{resu} 
\proof By definition of the set $z_{1} \cdot I_{1}+ \ldots + z_{r} \cdot I_{r}$ we
remark that:
\begin{equation}
  ^{t}\bz \cdot \boldsymbol{\a}=\{ \; ^{t}\bz\ba , \; \ba \in \boldsymbol{\a} \}.
\end{equation}
And the proof is thus obvious. In fact, it is a very simplified version of Beeck's Theorem (see e.g.
\cite{kn:Neumaier90} and references therein). \cqfd \\

Let us also recall the following rules of interval calculus (see
e.g. \cite{kn:Neumaier90}): \\

\noi
\textsf{AI1}. Two arithmetical expressions which are equivalent in real arithmetic are equivalent in
interval arithmetic when every variable occurs only once on each side.

\noi
\textsf{AI2}. If $f$ and $g$ are two arithmetical expressions of variables
$x_{1} \in I_{1}, \ldots, x_{n} \in I_{n}$ where $I_{1}, \ldots, I_{n}$ are given
intervals, which are equivalent in real
arithmetic then the inclusion $f(I_{1},\ldots, I_{n}) \subseteq g(I_{1}, \ldots, I_{n})$
holds if every variable $x_{i}$ occurs only once in $f$.

\subsection*{Acknowledgment}
Author would like to thank Odile Bellenguez, Gilles Chabert, Chams Lahlou,
Olivier Le Corre and James Ledoux
for helpful discussions.

\section{Main results}
\label{sec-MR}

Let us consider a matrix $\bA \in \textsf{Mat}(\real,m,n)$ satisfying
ASSUMPTIONS (A) and (B) (see subsection~\ref{sub-pb}), and a vector
$\grasb \in \real^{m}$.

As a direct consequence of ASSUMPTION (B) we can suppose that
$\bA=\left(\bar{c} \bA_{1} \\ \bA_{2}
\ear \right)$ with 
$\bA_{1}:=\left(\bar{c} \ba_{1,.} \\ \vdots \\ \ba_{m-n,.}\ear
\right) \in \textsf{Mat}(\real,m-n,n)$ and $\bA_{2}:=\left(\bar{c}
\ba_{m-n+1,.} \\ \vdots \\ \ba_{m,.}\ear \right) \in
\textsf{Mat}(\real,n,n)$ assumed to be invertible. Its inverse
is denoted $\bA_{2}^{-1}$. For all $i \in [m]$, $\ba_{i,.}$ denotes 
the $i$th row of matrix $\bA$.

Finally, let us recall that $\bA^{+}$ denotes the Moore-Penrose
inverse of matrix $\bA$ (see Definition~\ref{defMP}).

    \begin{lemm}
      \label{lem1}
We have the following two logical equivalences:
    \[
 \forall \bc, \;   \bA \bA^{+} \bc = \bc \Leftrightarrow (\bA \bA^{+}-\bI_{m}) \bc = \bzero \Leftrightarrow
    \bU \bc =\bzero
 \]
 with:

 \begin{equation}
   \label{defU}
 \bU := \left(\bar{cc} \bI_{m-n} & -\bA_{1} \bA_{2}^{-1} \\
 \bzero_{n,m-n} & \bzero_{n,n}
 \ear \right)
 \end{equation}
\end{lemm}

    \proof. The first equivalence is obvious. So, let us prove
$(\bA \bA^{+}-\bI_{m}) \bc = \bzero \Leftrightarrow
    \bU \bc =\bzero$. \\

    Using \cite{kn:Hung-Mar75} we develop $\bA^{+}$ as follows:

\[
\bA^{+}= \left(\bar{c} \bA_{1} \\ \bA_{2} \ear\right)^{+} = \left(\bar{cc} \bK^{+} \; ^{t}\bA_{1} &
\bK^{+} \; ^{t}\bA_{2} \ear\right),
\]
with $\bK:= \; ^{t}\bA_{1} \bA_{1} + \; ^{t}\bA_{2} \bA_{2}$.

Because $\bA_{2}$ is invertible, $^{t}\bA_{2} \bA_{2}$ is symmetric
invertible and thus $\bK$ is invertible. Hence, $\bK^{+}=\bK^{-1}$.

Now, we have by block multiplication of matrices:
\[
\bar{ll}
\bA \bA^{+} -\bI_{m} & = \left(\bar{c} \bA_{1} \\ \bA_{2}
\ear \right) \; \left(\bar{cc} \bK^{+} \; ^{t}\bA_{1} &
\bK^{+} \; ^{t}\bA_{2} \ear\right) - \left(\bar{cc} \bI_{m-n} & \bzero_{m-n,n} \\
\bzero_{n,m-n} & \bI_{n}
\ear \right) \\
\mbox{} & = \left(\bar{cc}  \bA_{1} \bK^{+} \; ^{t}\bA_{1} - \bI_{m-n} & \bA_{1} \bK^{+} \; ^{t}\bA_{2}\\
\bA_{2} \bK^{+} \; ^{t}\bA_{1} &  \bA_{2} \bK^{+} \; ^{t}\bA_{2} - \bI_{n}
\ear \right)
\ear
\]
Using block linear elimination we have (last row block multiplied by $\bA_{1} \bA_{2}^{-1}$ and
then substracting to first row block):
\[
\left(\bar{cc}  \bI_{m-n} & -\bA_{1} \bA_{2}^{-1}\\
\bA_{2} \bK^{+} \; ^{t}\bA_{1} &  \bA_{2} \bK^{+} \; ^{t}\bA_{2} - \bI_{n}
\ear \right)
\]
Left multiplying the last row block by $\bK \bA_{2}^{-1}$ one has:
\[
\left(\bar{cc}  \bI_{m-n} & -\bA_{1} \bA_{2}^{-1}\\
^{t}\bA_{1} &   ^{t}\bA_{2} - \bK \bA_{2}^{-1}
\ear \right)
\]
Right multypling the first row block by $-\; ^{t}\bA_{1}$ and adding to the second
row block we obtain:
\[
\left(\bar{cc}  \bI_{m-n} & -\bA_{1} \bA_{2}^{-1}\\
\bzero_{n,m-n}&   ^{t}\bA_{1} \bA_{1} \bA_{2}^{-1}  + \; ^{t}\bA_{2} - \bK \bA_{2}^{-1}
\ear \right).
\]
Now, we just have to note that:

\[
\bar{ll}
^{t}\bA_{1} \bA_{1} \bA_{2}^{-1}  + \; ^{t}\bA_{2} - \bK \bA_{2}^{-1} & = (^{t}\bA_{1} \bA_{1} +
\; ^{t}\bA_{2}\bA_{2}- \bK) \bA_{2}^{-1} \\
\mbox{} & = (\bK - \bK) \bA_{2}^{-1} \\
\mbox{} & = \bzero_{n,n}.
\ear
\]

Conversely, assume that $\bU \bc = \bzero$. This equality is equivalent to

\begin{equation}
   \label{systeme}
\bc_{1} = \bR \bc_{2},
  \end{equation}
with  $\bc_{1}:= \;^{t}(c_1, \ldots, c_{m-n})$, $\bc_{2}:=\;^{t}(c_{m-n+1}, \ldots, c_{m})$ and 
\begin{equation}
  \label{eqdefR}
\bR:= \bA_{1} \bA_{2}^{-1}.
  \end{equation}
Now, every $\bc$ such that $\bU \bc = \bzero$ has the form: $\hat{\bR} \bc_{2}$, with
$\hat{\bR}:= \left(\bar{c} \bR \\ \bI_{n} \ear \right)$. It remains to check that
$\bA \bA^{+} \hat{\bR}= \hat{\bR}$. We develop the computation as follows.

\[
\bar{ll}
\bA \bA^{+} \hat{\bR} & = \left(\bar{cc}  \bA_{1} \bK^{+} \; ^{t}\bA_{1}  & \bA_{1} \bK^{+} \; ^{t}\bA_{2}\\
\bA_{2} \bK^{+} \; ^{t}\bA_{1} &  \bA_{2} \bK^{+} \; ^{t}\bA_{2}
\ear \right) \; \left(\bar{c} \bR \\ \bI_{n} \ear \right) \\
\mbox{} & = \left(\bar{c}  \bA_{1} \bK^{+} \; ^{t}\bA_{1} \bR + \bA_{1} \bK^{+} \; ^{t}\bA_{2}\\
\bA_{2} \bK^{+} \; ^{t}\bA_{1} \bR + \bA_{2} \bK^{+} \; ^{t}\bA_{2} \ear \right).
\ear
\]

Now, we have:

\[
\bar{ll}
\bA_{1} \bK^{+} \; ^{t}\bA_{1} \bR + \bA_{1} \bK^{+} \; ^{t}\bA_{2} & = \bA_{1} \bK^{+} \; ^{t}\bA_{1} \bA_{1}
\bA_{2}^{-1} + \bA_{1} \bK^{+} \; ^{t}\bA_{2} \\
\mbox{ } &=\bA_{1} \bK^{+} \; (\; ^{t}\bA_{1} \bA_{1} + \; ^{t}\bA_{2} \bA_{2}) \bA_{2}^{-1} \\
\mbox{} & = \bA_{1} \bK^{+} \bK \bA_{2}^{-1} \\
\mbox{} &= \bA_{1} \bA_{2}^{-1}.
\ear
\]
And
\[
\bar{ll}
\bA_{2} \bK^{+} \; ^{t}\bA_{1} \bR + \bA_{2} \bK^{+} \; ^{t}\bA_{2} & = \bA_{2} \bK^{+} \; ^{t}\bA_{1} \bA_{1}
\bA_{2}^{-1} +  \bA_{2} \bK^{+} \; ^{t}\bA_{2} \\
\mbox{} &= \bA_{2} \bK^{+} (\; ^{t}\bA_{1} \bA_{1} + \; ^{t}\bA_{2} \bA_{2}) \bA_{2}^{-1} \\
\mbox{} &=  \bA_{2} \bK^{+} \bK \bA_{2}^{-1} \\
\mbox{} &= \bI_{n}.
\ear
\]
Thus, the result is proved. \cqfd \\

Let $^{t}\ba$ be any row vector of matrix $\bA_{1}$. Let $\tilde{\bA}_{2}$ be any submatix
of $\bA_{2}$ such that $\tilde{\bA}_{2} \in \textsf{Mat}(\real,k,n)$ for some $k \in [n]$. Let
$\tilde{\bA}:=\left(\bar{c} ^{t}\ba \\ \tilde{\bA}_{2} \ear \right)$. Finally, let us denote
$\tilde{\bc}:=\left(\bar{c} c \\ \tilde{\bc}_{2} \ear \right)$ with $c \in \real$ and
$\tilde{\bc}_{2} \in \real^{k}$.
Then

\begin{lemm}
  \label{lem2}

\[
\forall \tilde{\bc}, \;
(\tilde{\bA} \tilde{\bA}^{+} - \bI_{k+1}) \tilde{\bc} =\bzero \Leftrightarrow \tilde{\bU} \tilde{\bc} =\bzero,
\]
with
\begin{equation}
  \label{defUtilde}
\tilde{\bU}:=\left(\bar{cc} 1 & - \; ^{t}\ba \tilde{\bA}_{2}^{+} \\
\bzero_{k,1} & \bzero_{k,k}
\ear\right).
\end{equation}
  \end{lemm}
\proof. Let us denote $\bdelta:= \tilde{\bA}_{2}^{+} \ba$, and $\bgamma:= \ba - \tilde{\bA}_{2} \bdelta$.
Because $\tilde{\bA}_{2}$ is a $k \times n$-submatrix of the invertible
matrix $\bA_{2}$ the matrix $\tilde{\bA}_{2}^{+}$ is the right inverse of $\tilde{\bA}_{2}$, thus
$\tilde{\bA}_{2} \tilde{\bA}_{2}^{+}= \bI_{k}$. So, we are in the case where
$\bgamma = \bzero$. And we apply \cite[section 4]{kn:Greville60} to the block matrix $\tilde{\bA}$
to obtain:

\[
\tilde{\bA}^{+} = \left(\bar{cc} h^{-1} \tilde{\bA}_{2}^{+}\; ^{t}\tilde{\bA}_{2}^{+} \ba &
\tilde{\bA}_{2}^{+} - h^{-1} \tilde{\bA}_{2}^{+}\; ^{t}\tilde{\bA}_{2}^{+} \ba \; ^{t}\ba
\tilde{\bA}_{2}^{+} \ear \right),
\]
where $h:= 1 + \; ^{t}\bv \bv$ with $^{t}\bv:= \;^{t}\ba \tilde{\bA}_{2}^{+}$ and the
matrix $^{t}\tilde{\bA}_{2}^{+}$ is defined as the matrix $^{t}(\tilde{\bA}_{2}^{+})=
(^{t}\tilde{\bA}_{2})^{+}$. Now, we have:

\[
\tilde{\bA} \tilde{\bA}^{+} = \left(\bar{cc} h^{-1} \; ^{t}\bv \bv & ^{t}\bv - h^{-1} \; ^{t}\bv
^{t}\tilde{\bA}_{2}^{+} \ba \; ^{t}\bv \\
h^{-1} \tilde{\bA}_{2} \tilde{\bA}_{2}^{+}\bv  & \tilde{\bA}_{2} \tilde{\bA}_{2}^{+} - h^{-1}
\tilde{\bA}_{2} \tilde{\bA}_{2}^{+} \bv \; ^{t}\bv
\ear \right).
\]
After some manipulations one obtains that:
\[
\bar{ll}
\tilde{\bA} \tilde{\bA}^{+} -\bI_{k+1} &= \tilde{\bA} \tilde{\bA}^{+} - \left(\bar{cc} 1 & \bzero_{1,k} \\ \bzero_{k,1} & \bI_{k} \ear \right) \\
\mbox{} & = \left(\bar{cc} -h^{-1} & h^{-1} \; ^{t}\bv \\
h^{-1}\bv & -h^{-1} \bv \; ^{t}\bv
\ear \right).
\ear
\]
Multiplying by $-h$ the first row of the above matrix one has the following
matrix:
\[
\left(\bar{cc} 1 & -\; ^{t}\bv \\
h^{-1}\bv & -h^{-1} \bv \; ^{t}\bv
\ear \right).
\]
Then, left multiplying the first row by $-h^{-1} \bv$ and adding to the second row we have:
\[
\left(\bar{cc} 1 & -\; ^{t}\bv \\
\bzero_{k,1} & \bzero_{k,k}
\ear \right).
\]

Conversely, the system of equations $\tilde{\bU} \tilde{\bc} =\bzero$ is equivalent to
$c = \; ^{t}\bv \tilde{\bc}_{2} = \; ^{t}\ba \tilde{\bA}_{2}^{+} \tilde{\bc}_{2}$. Thus, we have to
prove that matrix
$\hat{\bR}:=\left(\bar{c} ^{t}\bv \\ \bI_{k} \ear \right)$ satisfies:
$\tilde{\bA} \tilde{\bA}^{+}\hat{\bR} = \hat{\bR}$.

We then, have:

\[
\bar{ll}
\tilde{\bA} \tilde{\bA}^{+}\hat{\bR} & = \left(\bar{cc} h^{-1} \; ^{t}\bv \bv & ^{t}\bv - h^{-1} \; ^{t}\bv
^{t}\tilde{\bA}_{2}^{+} \ba \; ^{t}\bv \\
h^{-1} \tilde{\bA}_{2} \tilde{\bA}_{2}^{+}\bv  & \tilde{\bA}_{2} \tilde{\bA}_{2}^{+} - h^{-1}
\tilde{\bA}_{2} \tilde{\bA}_{2}^{+} \bv \; ^{t}\bv
\ear \right) \; \left(\bar{c} ^{t}\bv \\ \bI_{k} \ear \right) \\
\mbox{} &= \left(\bar{c} h^{-1} \; ^{t}\bv \bv \; ^{t}\bv + ^{t}\bv - h^{-1} \; ^{t}\bv
^{t}\tilde{\bA}_{2}^{+} \ba \; ^{t}\bv \\
h^{-1} \tilde{\bA}_{2} \tilde{\bA}_{2}^{+}\bv \; ^{t}\bv + \tilde{\bA}_{2} \tilde{\bA}_{2}^{+} - h^{-1}
\tilde{\bA}_{2} \tilde{\bA}_{2}^{+} \bv \; ^{t}\bv
\ear \right).
\ear
\]
Noticing that $\bv=\; ^{t}\tilde{\bA}_{2}^{+} \ba$ we have:
\[
h^{-1} \; ^{t}\bv \bv \; ^{t}\bv + ^{t}\bv - h^{-1} \; ^{t}\bv
^{t}\tilde{\bA}_{2}^{+} \ba \; ^{t}\bv = h^{-1} \; ^{t}\bv \bv \; ^{t}\bv + ^{t}\bv - h^{-1} \; ^{t}\bv \bv \; ^{t}\bv
= \; ^{t}\bv.
\]
Noticing that $\tilde{\bA}_{2} \tilde{\bA}_{2}^{+}= \bI_{k}$ we have:
\[
h^{-1} \tilde{\bA}_{2} \tilde{\bA}_{2}^{+}\bv \; ^{t}\bv + \tilde{\bA}_{2} \tilde{\bA}_{2}^{+} - h^{-1}
\tilde{\bA}_{2} \tilde{\bA}_{2}^{+} \bv \; ^{t}\bv = h^{-1} \bv \; ^{t}\bv + \bI_{k} - h^{-1}
\bv \; ^{t}\bv = \bI_{k}.
\]
And the result is proved. \cqfd \\

For $i \in [m]$ we define:
\begin{equation}
  \label{eqdefLi}
\mathcal{L}_{i}:= \{\bx \in \real^{n}: \ba_{i, .} \bx \leq b_{i}\},
\end{equation}
recalling that $\ba_{i,.}$ denotes the $i$th row vector of matrix $\bA$.

Recall $\mathcal{C}_{i}=[-\infty, b_{i}]$, $i \in [m]$. And let us
define the vector of the intervals $\mathcal{C}_{i}$, $i \in [m]$ by
$\boldsymbol{\c}:= \; ^{t}(\mathcal{C}_{1}, \ldots, \mathcal{C}_{m})$. That is
\begin{equation}
  \label{defc}
\boldsymbol{\c} = \{\bc: \; \bminfty \leq \bc \leq \grasb \}.
\end{equation}

For $i \in [m-n]$ we define $B_{i}:= \{j \in [m]: u_{i,j} \neq 0 \}$ recalling
that $\bU =[u_{i,j}]$ is defined by (\ref{defU}).

    \begin{theo}
      \label{thm-main}

    For all $i \in [m-n]$ we have:

    \[
    0 \in \bu_{i,.} \cdot \boldsymbol{\c} \Leftrightarrow \cap_{j \in B_{i}} \mathcal{L}_{j} \neq \emptyset.
    \]
\end{theo}
\proof Without loss of generality we can assume that $i=m-n$. So, that:
\[
\bu_{m-n,.} =(\bzero_{1,m-n-1},1, - \ba_{m-n,.} \bA_{2}^{-1}).
\]
Which could be rewritten as:
\begin{equation}
  \label{defum-n}
\bu_{m-n,.} =(\bzero_{1,m-n-1},1, - \ba_{m-n,.} \bA_{2}^{-1} \tilde{\bPi}).
\end{equation}
Where $\tilde{\bPi}=[\tilde{\pi}_{i,j}] \in \textsf{Mat}(\real,n,n)$ such that:
$\tilde{\pi}_{i,j}=1$ if $i=j$ and $u_{m-n,i} \neq 0$, and $0$ otherwise.

Now, let us define $\bPi=[\pi_{i,j}] \in \textsf{Mat}(\real,m,m)$ by:

\[
\pi_{i,j} := \left\{\bar{ll} 0 & \mbox{ if $i \neq j$ or $i=j \leq m-n-1$} \\
1 & \mbox{ if $i=j= m-n$} \\
\tilde{\pi}_{i,j} & \mbox{ otherwise}.
\ear \right.
\]

By definition of $\bPi$, we remark that:

\[
\bar{ll}
\cap_{j \in B_{i}} \mathcal{L}_{j} \neq \emptyset & \Leftrightarrow \exists \bx, \bPi \bA \bx \leq \bPi \grasb \\
\mbox{} & \Leftrightarrow \exists \bx, \exists \bc \leq \grasb, \bPi \bA \bx = \bPi \bc \\
\mbox{} & \Leftrightarrow \exists \bc \leq \grasb, (\bPi \bA)(\bPi \bA)^{+} \bPi \bc = \bPi \bc \\
\mbox{} & \Leftrightarrow \exists \bc \leq \grasb, ((\bPi \bA)(\bPi \bA)^{+} \bPi - \bPi) \bc= \bzero.
\ear
\]
Let $k:=|B_{m-n}|$ be the number of elements of the set $B_{m-n}$. By renumbering the lines of matrix
$\bA$ we can assume that

\[
\bPi \bA = \left(\bar{c} \bzero_{m-n-1,n} \\ \ba_{m-n,.} \\ \tilde{\bPi} \bA_{2} \ear \right),
\]
with
\[
\tilde{\bPi}=\left(\bar{cc} \bI_{k} & \bzero_{k,n-k} \\
\bzero_{n-k,k} & \bzero_{n-k,n-k}
\ear \right).
\]

Because of the expression of vector $\bu_{m-n,.}$ (see (\ref{defum-n})) we
define $\bX \in \textsf{Mat}(\real,n,n)$
as $\bX:=\bA_{2}^{-1} \tilde{\bPi}$. It is easy to check that $\bX$ is a
$\{1,2\}$-inverse (see Definition~\ref{defMP}) of the matrix $\tilde{\bPi} \bA_{2}$ 
such that: $(\tilde{\bPi} \bA_{2}) \bX = \tilde{\bPi}$. The latter equation
means that because $\tilde{\bPi} \bA_{2} = \left(\bar{c} \tilde{\bA}_{2} \\ \bzero_{n-k,n} \ear\right)$ where $\tilde{\bA}_{2}$ is a submatrix of $\bA_{2}$,
we have $\bX =\left(\bar{cc} \tilde{\bA}_{2}^{+} & \bzero_{k,n-k} \ear \right)$
where $\tilde{\bA}_{2}^{+}$ is the right inverse of $\tilde{\bA}_{2}$.

Noticing matrix $\bPi \bA$ has the form $\left(\bar{c} \bzero_{m-n-1,n} \\ \bY \\
\bzero_{n-k,n} \ear\right)$
with $\bY:=\left(\bar{c} \ba_{m-n,.} \\ \tilde{\bA}_{2} \ear \right)$, its Moore-Penrose
inverse is then the matrix $\left(\bar{ccc} \bzero_{n,m-n-1} & \bY^{+} & \bzero_{n,n-k} \ear \right)$. Hence,
\[
\bPi \bA (\bPi \bA)^{+} = \left(\bar{ccc} \bzero_{m-n-1,m-n-1} & \bzero_{m-n-1,k+1} & \bzero_{m-n-1,n-k} \\
\bzero_{k+1,m-n-1} & \bY \bY^{+} & \bzero_{k+1,n-k} \\
\bzero_{n-k,m-n-1} & \bzero_{n-k,k+1} & \bzero_{n-k,n-k}
\ear \right).
\]
By definition of matrix $\bPi$ we then have:
\[
\bPi \bA (\bPi \bA)^{+} \bPi -\bPi = \left(\bar{ccc} \bzero_{m-n-1,m-n-1} & \bzero_{m-n-1,k+1} & \bzero_{m-n-1,n-k} \\
\bzero_{k+1,m-n-1} & \bY \bY^{+}-\bI_{k+1} & \bzero_{k+1,n-k} \\
\bzero_{n-k,m-n-1} & \bzero_{n-k,k+1} & \bzero_{n-k,n-k}
\ear \right).
\]
We can focus
our attention on matrix $\bY=\left(\bar{c} \ba_{m-n,.} \\ \tilde{\bA}_{2} \ear \right)$ and
apply Lemma~\ref{lem2} to obtain that the system of equations 
$(\bPi \bA (\bPi \bA)^{+} \bPi -\bPi) \bc = \bzero$ is equivalent to $\bu_{m-n, .} \bc=0$.
Hence the result is now proved. \cqfd \\

\begin{theo}
  \label{thm-main2}
  The polyhedron $\mathcal{P}(\bA, \grasb)$ is not empty iff 

  \[
\mathfrak{B}: \; \forall \bk, 0 \in \; ^{t}\bk \bU \cdot \boldsymbol{\c}.
  \]

  \end{theo}
\proof Using the characterization of non emptyness of polyhedron $\mathcal{P}(\bA, \grasb)$
(see (\ref{eq-nonempty})) and Lemma~\ref{lem1} we have:
\[
\mathcal{P}(\bA, \grasb) \neq \emptyset \Leftrightarrow \mathfrak{A}: \exists \bc \leq \grasb, \bU \bc =\bzero,
\]
recalling that $\bU = \left(\bar{cc} \bI_{m-n} & -\bR \\
 \bzero_{n,m-n} & \bzero_{n,n}
 \ear \right)$ with $\bR=\bA_{1} \bA_{2}^{-1}$. Thus,
 we have to prove: $\mathfrak{A} \Leftrightarrow \mathfrak{B}$. \\

 \noi
By application of Result~\ref{res-Beeck} $\mathfrak{A} \Rightarrow \mathfrak{B}$. \\

 \noi
 Let us prove $\mathfrak{B} \Rightarrow \mathfrak{A}$ by absurd. Thus, assume $\mathfrak{B}$
 and $\overline{\mathfrak{A}}: \;
 \forall \bc \leq \grasb$, $\exists i, \bu_{i,.} \bc \neq 0$.

 Let $\bc \leq \grasb$ then by $\overline{\mathfrak{A}}$ there exists
 $i \in [m]$ such that $\bu_{i,.} \bc \neq 0$. But, take $\bk=\be^{i}$,
 then by $\mathfrak{B}$ we have: $0 \in \; ^{t}\be^{i} \bU \cdot \boldsymbol{\c}=\bu_{i,.} \cdot \boldsymbol{\c}$ which means that $\exists \bc^{1} \leq \grasb$
 such that $\bu_{i,.} \bc^{1}=0$ (see Result~\ref{res-Beeck}). Let us
 define $\mathcal{C}:=[-\infty, b_{1}] \times \cdots \times [-\infty, b_{m}]$.
 
 The vector $\bc^{1}$ is $\leq \grasb$ and such that $\exists i_{1} \in [m]$ with
 $\bu_{i_{1},.} \bc^{1} \neq 0$ (by $\overline{\mathfrak{A}}$). But by
 $\mathfrak{B}$ there exists $\bc^{2}$ such that
 $\bu_{i_{1},.} \bc^{2} =0$. So, we construct a $[m] \times \mathcal{C}$-valued
 series $\sigma:= \{(i_{0}=i, ^{t}\bc^{0}= \;^{t}\bc), (i_{1}, ^{t}\bc^{1}), (i_{2}, ^{t}\bc^{2}), (i_{3}, ^{t}\bc^{3}), \ldots \}$
 which has the following property $\forall n \geq 0$: \\

 \noi
 (p). $\bu_{i_{n}, .} \bc^{n} \neq 0$ and $\bu_{i_{n},.} \bc^{n+1}=0$. \\

 The set $[m] \times \mathcal{C}$ is the cartesian product of compact
 spaces thus it is a compact space. Hence the sequence $\sigma$ admits
 a subsequence $\varphi.\sigma:=\{(i_{\varphi(k)}, \bc^{\varphi(k)}),
 k \in \integer\}$ with $\varphi: \integer \rightarrow \integer$
 strictly increasing and such that $\lim_{k \rightarrow \infty}
 \varphi(k)=\infty$. And the subsequence $\varphi.\sigma$ admits a
 limit $(\iota, \boldsymbol{\ell}):= \lim_{k \rightarrow \infty}
 (i_{\varphi(k)}, \bc^{\varphi(k)})$. The point
 $(\iota, \boldsymbol{\ell})$ is an element of $[m] \times
 \mathcal{C}$ which must satisty property (p) by construction. Thus,
 we obtain a contradiction. And the result is now proved. \cqfd \\

  \begin{rem}
   \label{rem2}
   Due to the structure of the matrix $\bU$
   we can restrict our attention to all $m-n$ dimensional vectors
  $\bk'$ such that
  \begin{equation}
    \label{eq-pour-algo}
0 \in \; ^{t}\bk' \bG \cdot \boldsymbol{\c},
  \end{equation}
  where matrix $\bG$ is:
  \begin{equation}
\bG := \left(\bar{cc} \bI_{m-n} & -\bR \ear \right).
    \end{equation}
  \end{rem}

The main problem is then to enumerate only the relevant vectors $\bk'$ such
that $ 0 \in \; ^{t}\bk' \bG$. The following result adresses
this problem.

Let $\mathcal{B}:=\{\be^{i}\}$ be the canonical basis of $\real^{m-n}$.
Let us define
\begin{equation}
  \label{k'jii'}
  \bk'(j,i,i'):=  -r_{i',j} \be^{i} + r_{i,j} \be^{i'},
\end{equation}
for all $j \in [n]$, $i \in [m-n-1]$ and $i'=i+1$ to $m-n$.

Let us define
\begin{equation}
  \label{kerR}
 \textsf{ker}^{l}(\bR):=\{\bk': \;  ^{t}\bk' \bR =\; ^{t}\bzero \}.
\end{equation}
The set $\mathcal{K}$ denotes a basis of $\textsf{ker}^{l}(\bR)$ if
$\textsf{ker}^{l}(\bR) \neq \{\bzero \}$ and $\{\bzero \}$
otherwise. 

Let us define $\grasb_{1}= \left(\bar{c} b_{1} \\ \vdots \\ b_{m-n}\ear \right)$ and
$\grasb_{2}= \left(\bar{c} b_{m-n+1} \\ \vdots \\ b_{m}\ear \right)$. Let us denote
$(\grasb_{1})^{\perp}$ a basis of the set $\{\bk': \;^{t}\bk' \grasb_{1}=0\}$. And let us
denote $(\bR \grasb_{2})^{\perp}$ a basis of the set $\{\bk': \; ^{t}\bk' (\bR \grasb_{2})=0\}$.

\begin{theo}
  \label{th-main-3}
If for all $i \in [m-n]$:

\[
0 \in \;^{t}\be^{i} \bG \cdot \boldsymbol{\c},
\]
and
\[
\forall \bk' \in \mathcal{K}, \; 0 \in \; ^{t}\bk'
\bG \cdot \boldsymbol{\c},
\]
and
\[
\forall \bk' \in (\grasb_{1})^{\perp}, \; 0 \in \; ^{t}\bk'
\bG \cdot \boldsymbol{\c},
\]
and
\[
\forall \bk' \in (\bR\grasb_{2})^{\perp}, \; 0 \in \; ^{t}\bk'
\bG \cdot \boldsymbol{\c},
\]
%with
%\begin{equation}
% \textsf{ker}^{l}(\bR):=\{\bk': \;  ^{t}\bk' \bR =\; ^{t}\bzero \}.
%\end{equation}
and for all $j \in [n]$, $i \in [m-n-1]$ and $i'=i+1$ to $m-n$:
  \[
0 \in \; ^{t}\bk'(j,i,i') \bG \cdot \boldsymbol{\c}.
\]
%and for all $i \in [m-n-1]$ and $i'=i+1$ to $m-n$ such that
%$\mathcal{J}(i,i') = \emptyset$:
%\[
%0 \in \; ^{t}\bk'(i,i') \bG \cdot \boldsymbol{\c}.
%\]
\noi
Then
\[
\forall \bk', 0 \in \; ^{t}\bk' \bG \cdot \boldsymbol{\c}.
\]
\end{theo}
\proof First, let us remark that by definition of $ \boldsymbol{\c}$
(see (\ref{defc})) we have:

\[
^{t}\bk' \bG \cdot \boldsymbol{\c}= \left\{\bar{ll} \mbox{$[-\infty, \; ^{t}\bk' \bG \grasb]$} & \mbox{if $^{t}\bk' \bG \geq \; ^{t}\bzero$} \\
\mbox{$[\; ^{t}\bk' \bG \grasb, +\infty]$} & \mbox{if $^{t}\bk' \bG \leq \; ^{t}\bzero$} \\
\mbox{$[-\infty, +\infty]$} & \mbox{otherwise}.
\ear \right.
\]
Noticing that $^{t}\bk' \bG \leq \; ^{t}\bzero \Leftrightarrow -\;^{t}\bk' \bG
\geq \; ^{t}\bzero$ we can focus our attention on the cone:
\[
\mathcal{G}:=\{\bk': \;^{t}\bk' \bG
\geq \; ^{t}\bzero \}.
\]
If $\mathcal{G}=\{\bzero\}$ then the result is obviously true. Thus, let us
assume that the cone $\mathcal{G} \neq \{\bzero\}$. In this case we have:

\[
\forall \bk' \in \mathcal{G}, 0 \in \; ^{t}\bk' \bG \cdot \boldsymbol{\c} \Leftrightarrow
0 \in \cap_{\bk' \in \mathcal{G}}[-\infty, \;^{t}\bk' \bG \grasb] \Leftrightarrow
\forall \bk' \in \mathcal{G}, \; ^{t}\bk' \bG \grasb \geq 0.
\]
The last equivalence is then equivalent to
$\min_{\bk' \in \mathcal{G}}(\;^{t}\bk' \bG \grasb) \geq 0$.

Now, let us remark that:
\[
f(\bk'):= \;^{t}\bk' \bG \grasb = \; ^{t}\bk' \grasb_{1} - ^{t}\bk' \bR \grasb_{2}.
\]
Due to the structure
of matrix $\bG$ the inequality $^{t}\bk' \bG \geq \; ^{t}\bzero$ implies
$\bk' \geq \bzero$. The different possibilities for the choice of vector
$\bk' \in \mathcal{G}$ in the function $f$ are as follows:

\ben
\item $\bk'=\be^{i}$, $i \in [m-n]$. In such case $f=\bg_{i,.} \grasb$.
\item $\bk' \in \mathcal{K}$ or $\bk' \in (\bR \grasb_{2})^{\perp}$. In this case
  $f= \;^{t}\bk' \grasb_{1}$.
\item $\bk' \in (\grasb_{1})^{\perp}$. And then, $f= \; - ^{t}\bk' \bR \grasb_{2}$.
\item Finally, we can eliminate variables $b_{m-n+j}$, $j \in [n]$, between the rows
  $i$ and $i'$ of matrix $\bR$ for $i \in [m-n-1]$, $i'=i+1, \ldots, m-n$. In this
  case $f= \; ^{t}\bk'(j,i,i') \bG \grasb$.
\een

Hence, the result. \cqfd \\

Based on the previous results definitions and notations we provide the
following algorithm for testing if a polyhedron is empty or not.

\section*{Algorithm}

\bit
\item Inputs: $\bA \in \textsf{Mat}(\real,m,n)$ satisfying ASSUMPTIONS (A) and (B), and
  vector $\grasb \in \real^{m}$
\item Output: answer to the question `` is
  $\mathcal{P}(\bA, \grasb):=\{ \bx \in \real^{n}: \bA x \leq \grasb\}$ empty ?''

\item[0.] Put matrix $\bA$ in the form $\left(\bar{c} \bA_{1} \\ \bA_{2}\ear \right)$ with
  $\bA_{2}$ invertible.

\item[1.]  Compute $\bA_{2}^{-1}$.

\item[2.] Compute $\bG := \left(\bar{cc} \bI_{m-n} & -\bR \ear \right)$, with
  $\bR:= \bA_{1} \bA_{2}^{-1}$.

\item[3.] Compute $(\grasb_{1})^{\perp}$ a basis of the set $\{\bk': \;^{t}\bk' \grasb_{1}=0\}$. 
  \[
  (\grasb_{1})^{\perp}_{+}:= ((\grasb_{1})^{\perp} \cup (-(\grasb_{1})^{\perp})) \cap \mathcal{G}.
  \]
\item[4.] Compute $(\bR \grasb_{2})^{\perp}$ a basis of the set $\{\bk': \; ^{t}\bk' (\bR \grasb_{2})=0\}$.
  \[
(\bR\grasb_{2})^{\perp}_{+}:= ((\bR\grasb_{2})^{\perp} \cup (-(\bR\grasb_{2})^{\perp})) \cap \mathcal{G}.
  \]
  
\item[5.] Compute $\mathcal{K}$ a basis of $\textsf{ker}^{l}(\bR)$.
  \[
\mathcal{K}_{+}:= (\mathcal{K} \cup (-\mathcal{K})) \cap \mathcal{G}.
  \]
  
\item[6.] For $\bk' \in (\grasb_{1})^{\perp}_{+} \cup (\bR \grasb_{2})^{\perp}_{+} \cup \mathcal{K}_{+}$
  \bit
\item if $0 \notin \; ^{t}\bk' \bG \cdot \boldsymbol{\c}$ exit: polyhedron $\mathcal{P}(\bA, \grasb) = \emptyset$
  \eit
\item[] EndFor
  
\item[7.] For $i=1$ to $m-n$
  \bit
\item if $0 \notin \; ^{t}\be^{i} \bG \cdot \boldsymbol{\c}$ exit: polyhedron $\mathcal{P}(\bA, \grasb) = \emptyset$
  \eit
\item[] EndFor
  
\item[8.] For $j=1$ to $n$

  \bit
\item For $i=1$ to $m-n$
  \bit
\item For $i'=i+1$ to $m-n$
  \bit
\item $\bk' := -r_{i',j} \be^{i} + r_{i,j} \be^{i'}$
\item  if $0 \notin \; ^{t}\bk' \bG \cdot \boldsymbol{\c}$ exit: polyhedron $\mathcal{P}(\bA, \grasb) = \emptyset$
  \eit
\item EndFor
  \eit
  \item EndFor
  \eit

\item[] EndFor
  \eit

\section{Discussion about ASSUMPTIONS (A) and (B)}
\label{sec-Discuss}

\subsection{ASSUMPTION (A)}
\label{assA}

Independently of the definitions of a polyhedron (see subsection~\ref{assB}), if there
exists a row $\ba_{i,.}= \; ^{t}\bzero$ then the set $\mathcal{L}_{i}$ (see
(\ref{eqdefLi})) is thus defined as
\[
\mathcal{L}_{i}= \{\bx \in \real^{n}: \; ^{t}\bzero \bx \leq b_{i}\}.
\]
Then if $b_{i} <0$ $\mathcal{L}_{i} = \emptyset$ and thus the polyhedron
$\mathcal{P}(\bA, \grasb)= \emptyset$. Otherwise the inequality
$^{t}\bzero \bx \leq b_{i}$ can be removed.

\subsection{ASSUMPTION (B)}
\label{assB}
In this subsection we discuss several definitions of polyhedron which 
appear in linear programming. We consider a matrix
$\tilde{\bA} \in \textsf{Mat}(\real, \tilde{m}, \tilde{n})$ whose rank is
$\textsf{rk}(\tilde{\bA})=r$. We assume that $r \neq \tilde{n}$. And we
consider a vector $\tilde{\grasb} \in \real^{\tilde{m}}$. 

\bit

\item If the polyhedron is defined as the following set:
\[
\mathcal{P}(\tilde{\bA}, \tilde{\grasb}) := \{\bx \in \real^{\tilde{n}}:
\tilde{\bA} \bx \leq \tilde{\grasb}, \bx \geq \bzero \},
\]
then, we have:
\[
\tilde{\bA} \bx \leq \tilde{\grasb}, \bx \geq \bzero \Leftrightarrow 
\bA \bx \leq \grasb,
\]
with: $\bA:=\left(\bar{c}\tilde{\bA} \\ -\bI_{\tilde{n}} \ear\right) 
\in \textsf{Mat}(\real, \tilde{m}+\tilde{n}, \tilde{n})$ and 
$\grasb:=\left(\bar{c}\tilde{\grasb} \\ \bzero \ear\right) 
\in \real^{\tilde{m}+\tilde{n}}$. It is clear that 
$m:=\tilde{m}+\tilde{n}> n:=\tilde{n}$ and that $\bA$ has 
full column-rank (equal to $n$).

\item If the polyhedron is defined as the following set:
\[
\mathcal{P}(\tilde{\bA}, \tilde{\grasb}) := \{\bx \in \real^{\tilde{n}}:
\tilde{\bA} \bx = \tilde{\grasb}, \bx \geq \bzero \},
\]
then, we have:
\[
\tilde{\bA} \bx = \tilde{\grasb}, \bx \geq \bzero \Leftrightarrow 
\bA \bx \leq \grasb,
\]

with: $\bA:=\left(\bar{c}\tilde{\bA}\\ -\tilde{\bA} \\ -\bI_{\tilde{n}} 
\ear\right) \in \textsf{Mat}(\real, 2\; \tilde{m}+\tilde{n}, \tilde{n})$ 
and $\grasb:= \left(\bar{c}\tilde{\grasb}\\-\tilde{\grasb}
 \\ \bzero \ear\right) \in \real^{ 2 \; \tilde{m}+\tilde{n}}$. 
Once again, it is clear that 
$m:= 2 \; \tilde{m}+\tilde{n}> n:=\tilde{n}$ and that $\bA$ has 
full column-rank (equal to $n$).

\item Finally, if a polyhedron is defined as the set 
\[
\mathcal{P}(\tilde{\bA}, \tilde{\grasb}) := \{\bx \in \real^{\tilde{n}}:
\tilde{\bA} \bx \leq \tilde{\grasb} \}.
\]
Writing every element $\bx$ of 
$\real^{\tilde{n}}$ as: $\bx = \bx_{+} - \bx_{-}$ with $\bx_{+}, \bx_{-} 
\geq \bzero$ one has:
\[
\bar{ll}
\tilde{\bA} \bx \leq \tilde{\grasb} & \Leftrightarrow \left\{
\bar{ll} \tilde{\bA} (\bx_{+} - \bx_{-}) & \leq \tilde{\grasb} \\
\bx_{+}, \bx_{-} & \geq \bzero
\ear \right. \\
\mbox{} & \Leftrightarrow \left(\bar{cc} \tilde{\bA} & -\tilde{\bA} \\
                                         -\bI_{\tilde{n}} & \bzero_{\tilde{n},\tilde{n}} \\
                                             \bzero_{\tilde{n},\tilde{n}}   & -\bI_{\tilde{n}}
\ear \right) \; \left(\bar{c} \bx_{+} \\ \bx_{-} \ear \right) 
\leq \left(\bar{c} \tilde{\grasb} \\ \bzero \\ \bzero \ear \right)
\ear
\]

Let $m:= \tilde{m} + 2 \; \tilde{n}$ and $n:= 2\; \tilde{n}$. Then, the 
matrix $\bA:= \left(\bar{cc} \tilde{\bA} & -\tilde{\bA} \\
                                         -\bI_{\tilde{n}} & \bzero_{\tilde{n},\tilde{n}} \\
                                             \bzero_{\tilde{n},\tilde{n}}   & -\bI_{\tilde{n}}
\ear \right)$ is an element of $\textsf{Mat}(\real,m,n)$ such that 
its submatrix $\bA_{2}:=\left(\bar{cc} -\bI_{\tilde{n}} & \bzero_{\tilde{n},\tilde{n}} \\
                                             \bzero_{\tilde{n},\tilde{n}}   & -\bI_{\tilde{n}}
\ear \right)$ is clearly invertible. And thus, $\bA$ has full column rank 
$n$ with $m > n$.

\eit

\bibliographystyle{plain}
\bibliography{ref_lt}

\end{document}